\documentclass[12pt]{article}
\usepackage{amssymb}
\usepackage{amsbsy}
\begin{document}

\centerline{\Large\bf Mixed-Mean Inequality for Submatrix}



\vskip 20pt

\begin{center}

{  Lin Si, Suyun Zhao  }

\vskip 16pt

\begin{minipage}{11cm}
{\bf Abstract.} {\small For a $m\times n$ matrix $B=(b_{ij})_{m\times n}$ with
nonnegative entries $b_{ij}$ and any $k\times l-$submatrix $B_{ij}$ of $B$,
let $a_{B_{ij}}$ and $g_{B_{ij}}$ denote the arithmetic mean and
geometric mean of elements of $B_{ij}$ respectively. It is proved
that if $k$ is an integer in $(\frac{m}{2}, m]$ and $l$ is an
integer in $(\frac{n}{2}, n]$ respectively, then
$$\Big(\prod_{i=k,j=l\atop B_{ij}\subset B}a_{B_{ij}}\Big)^{\frac{1}{C_m^k\cdot C_n^l}}
\geq\frac{1}{C_m^k\cdot C_n^l}\Big(\sum_{i=k,j=l\atop
B_{ij}\subset B}g_{B_{ij}}\Big),$$ with equality if and only if
$b_{ij}$ is a constant for every $i,j$. }

\end{minipage}

\end{center}

\vskip 20pt

\centerline{{\large \bf 1.~~Introduction }}

\vskip 10pt

Let $x_1,...,x_n$ be positive real numbers, then the arithmetic-geometric mean inequality is
$$\frac{x_1+\cdot \cdot \cdot +x_n}{n}\geq\sqrt[n]{x_1\cdot \cdot \cdot x_n},$$
with equality if and only if $x_1=\cdot \cdot \cdot =x_n$.

There are many research articles devoted to the classical arithmetic-geometric mean inequality and its
generalizations([1],[2],[4],[9],[10]). The mixed-type arithmetic-geometric mean inequalities([3],[5],[6],[7],[8],[11]) are one of the most important branch in these generalizations.

In [6], Kedlaya established the
following mixed mean inequality, which was conjectured by F.Holland and was given an inductive proof by T.Matsuda[8]:

{\it The arithmetic mean of the numbers
$$x_1,\sqrt{x_1x_2},\sqrt[3]{x_1x_2x_3},...,\sqrt[n]{x_1x_2\cdot \cdot \cdot x_n}$$
does not exceed the geometric mean of the numbers
$$x_1, \frac{x_1+x_2}{2}, \frac{x_1+x_2+x_3}{3}, \frac{x_1+x_2+\cdot \cdot \cdot +x_n}{n},$$
with equality if and only if $x_1=\cdot \cdot \cdot =x_n$. }

In [3], Carlson established the
following mixed mean inequality:

{\it Let the arithmetic and geometric means of the real
nonnegative numbers $x_1, ..., x_n$ taken $n-1$ at a time be
denoted by
$$a_i=\frac{x_1+...+x_n-x_i}{n-1},
~g_i=\Big(\frac{x_1\cdot\cdot\cdot
x_n}{x_i}\Big)^{\frac{1}{n-1}}.$$ Then for $n\geq 3,$
$$\Big(a_1\cdot\cdot\cdot a_n\Big)^{\frac{1}{n}}\geq
\frac{g_1+...+g_n}{n},$$
with equality if and only if $x_1=\cdot \cdot \cdot =x_n$.}

In [7], Leng, Si and Zhu generalized the above result to any subsets
and established another mixed arithmetic-geometric mean inequality:

{\it Let $X=\{x_1,...,x_n~|~x_i>0, ~i=1,2,...,n\}.$ For $A\subset
X,$ let $a_A$ and $g_A$ denote the arithmetic mean and geometric
mean of all elements of $A$, respectively.  If $ k$ is an integer
in $(\frac{n}{2},n]$, then
$$(\prod_{|A|=k\atop A\subset X}a_{A})^{\frac{1}{C_n^k}}
\geq\frac{1}{C_n^k}(\sum_{|A|=k\atop A\subset X}g_{A}),\eqno(1)$$
with equality if and only if $x_1=...=x_n$}.

In this note, we established a new mixed arithmetic-geometric mean
inequality for submatix, which was an extension of the Carlson
inequality and also an extension of (1).

Our main result is the following theorem. \vskip 5pt

\noindent{\bf Theorem}~ {\it \small For a $m\times n$ matrix
$B=(b_{ij})_{m\times n}$ with
nonnegative entries $b_{ij}\geq 0$ and any $k\times
l-$submatrix $B_{ij}$ of $B$, let $a_{B_{ij}}$ and $g_{B_{ij}}$
denote the arithmetic mean and geometric mean of elements of
$B_{ij}$ respectively. If $k$ is an integer in $(\frac{m}{2}, m]$
and $l$ is an integer in $(\frac{n}{2}, n]$ respectively, then
$$\Big(\prod_{i=k,j=l\atop B_{ij}\subset B}a_{B_{ij}}\Big)^{\frac{1}{C_m^k\cdot C_n^l}}
\geq\frac{1}{C_m^k\cdot C_n^l}\Big(\sum_{i=k,j=l\atop
B_{ij}\subset B}g_{B_{ij}}\Big),\eqno(2)$$ with equality if and only if
$b_{ij}$ is a constant for every $i,j$}.

\vskip 8pt

\noindent{\it Remark} 1.~~~If $k\leq [\frac{m}{2}],$ for the matrix
$B=(b_{ij})_{m\times n}$, taking
$b_{11}=b_{12}=...=b_{1n}=b_{21}=b_{22}=...=b_{2n}=...=b_{k1}=b_{k2}=...=b_{kn}=1, ~b_{ij}=0$ for $k+1\leq i\leq m$, then the right-hand
side of (2) equals $1/({C_m^k\cdot C_n^l})$, but the left-hand side is zero. If $l\leq [\frac{n}{2}]$, by the same argument, one can get a contradiction. Hence the statement of Theorem fails for $k\leq [\frac{m}{2}]$ or $l\leq [\frac{n}{2}]$.

\vskip 5pt

\noindent{\it Remark} 2.~~Taking $k=m,l=n$ in Theorem,
inequality~(2) is just the classical arithmetic-geometric mean inequality.
For $m=1$ or $n=1$ in Theorem, inequality~(2) is just the inequality (1).

\vskip 5pt

\noindent{\it Remark} 3.~~The condition of our theorem is weaker than the ones of (1), because the infimum of $k\times l$ is $\frac{m}{2}\times \frac{n}{2}$, which is less than $\frac{1}{2}\times m\times n$, the half of the element number of the matrix
$B=(b_{ij})_{m\times n}$.

\vskip 20pt

\centerline{{\large \bf 2.~~Proof of Main Results }}

\vskip 10pt

Let $X$ denote the finite set with positive real numbers $x_1,x_2,...,x_n$ and let $X_i$ denote the subset of $X$ with $k$ elements
$x_{i_1},...,x_{i_k}$. The $r$ power mean of the elements of $X_i$ is denoted by $$m_{r}{(X_i)}=[\frac{1}{k}(x_{i_1}^r+...+x_{i_k}^r)]^{\frac{1}{r}},\eqno(3)$$
where $r>0$ is some real number. If $r=1$ in the above equality, we get $a_{X_i}$, the arithmetic mean of the elements of $X_i$. Let $r\rightarrow 0$ in (3), we get $g_{X_i}$, the geometric mean of the elements of $X_i$.

We first established the following lemma.

\noindent{\bf Lemma}~ {\it For any $m\times n$ matrix
$X=(x_{ij})_{m\times n}$ with
nonnegative entries $x_{ij}\geq 0$. Denote by $X_1,...,X_{C_m^k\cdot C_n^l}$ its all $k\times
l-$submatrix $X_{ij}$ of $X$. If $k$ is an integer in $(\frac{m}{2}, m]$
and $l$ is an integer in $(\frac{n}{2}, n]$ respectively, then
$$m_{r}{(X_i)}=\Big[\frac{1}{C_m^k\cdot C_n^l}\big[(m_r{(X_i\cap X_1)})^r+(m_r{(X_i\cap X_2)})^r+...
+(m_r{(X_i\cap X_{C_m^k\cdot C_n^l}}))^r\big]\Big]^{\frac{1}{r}}.\eqno(4)$$}
Here the intersection of $X_i\cap X_j, j=1,...,C_m^k\cdot C_n^l$, is just the general set intersection.
\vskip 5pt
\noindent{\it Proof.} It suffices to prove the equality
$$(m_{r}{(X_i)})^r=\frac{1}{C_m^k\cdot C_n^l}\big[(m_r{(X_i\cap X_1)})^r+(m_r{(X_i\cap X_2)})^r+...
+(m_r{(X_i\cap X_{C_m^k\cdot C_n^l}}))^r\big].\eqno(5)$$

Assume that $X_i=\{x_{i_1},...,x_{i_{k\times l}}\}$, then the left hand of (5) is that $$(m_{r}{(X_i)})^r=\frac{1}{k\times l}(x_{i_1}^r+...+x_{i_{k\times l}}^r).$$

Now we need show that the right hand of (5) is also the mean of $\{x_{i_1}^r,...,x_{i_{k\times l}}^r\}$ with the same coefficient $\frac{1}{k\times l}$.

Assume that the right hand of (5) is $c_{i_1}x_{i_1}^r+...+c_{i_{k\times l}}x_{i_{k\times l}}^r$, by the arbitrariness of $X_i$, we get that $c_{i_1}=...=c_{i_{k\times l}}$.

Since $k>\frac{m}{2}$ and $l>\frac{n}{2}$, we have
$$X_i\cap X_j\neq\Phi, ~j=1,2,...,C_m^k\cdot C_n^l.$$ Then the sum of all coefficients in $(m_r{(X_i\cap X_j)})^r$ is $1$, $~j=1,2,...,C_m^k\cdot C_n^l$.
As a result, the sum of all coefficients in the right hand of (5) is
$$\frac{1}{C_m^k\cdot C_n^l}(\underbrace{1+1+\cdot\cdot\cdot +1}_{C_m^k\cdot C_n^l})=1,$$
i.e., $c_{i_1}+\cdot\cdot\cdot+c_{i_{k\times l}}=1$.
Therefore $c_{i_1}=\cdot\cdot\cdot=c_{i_{k\times l}}=\frac{1}{k\times l}$. ~~~~~~~~~~~~~~~~~~~~~~~~~~~~~~~~~~~~$\square$

\noindent{\it Remark} 1.~~(4) is the generalization of related result in [7]. But the condition in Lemma is weaker, because the infimum of $k\times l$ is $\frac{m}{2}\times \frac{n}{2}$, which is less than $\frac{1}{2}\times m\times n$, the half of the element number of the matrix
$X=(x_{ij})_{m\times n}$.

\noindent{\it Remark} 2.~~In (4), if $r=1$, we have
$$a_{X_i}=\frac{1}{C_m^k\cdot C_n^l}(a_{X_i\cap X_1}+a_{X_i\cap X_2}+...
+a_{X_i\cap X_{C_m^k\cdot C_n^l}}).\eqno(6)$$
In (4), if $r\rightarrow 0$, we have
$$g_{X_i}=(g_{X_i\cap X_1}\cdot g_{X_i\cap X_2}\cdot ...
\cdot g_{X_i\cap X_{C_m^k\cdot C_n^l}})^{\frac{1}{C_m^k\cdot C_n^l}}.\eqno(7)$$

\vskip 20 pt
\noindent{\bf Proof of Theorem.}~

For a $m\times n$ matrix
$B=(b_{ij})_{m\times n}$ with
nonnegative entries $b_{ij}\geq 0$ and any $k\times
l-$submatrix $B_{ij}$ of $B$, which is composed by the lines $(i_1,i_2,...,i_k)$ of $B$ and by the rows $(j_1,j_2,...,j_l)$ of $B$, by the partial order of $(i_1,i_2,...,i_k)$ and $(j_1,j_2,...,j_l)$, we get
a subset chain $B_1,B_2,...,B_{C_m^k\cdot C_n^l}$

In Lemma, let $X$ be $B$ and $X_i$,$~j=1,2,...,C_m^k\cdot C_n^l$, be corresponding to $B_1,B_2,...,B_{C_m^k\cdot C_n^l}$, then (6) is
$$a_{B_i}=\frac{1}{C_m^k\cdot C_n^l}(a_{B_i\cap B_1}+a_{B_i\cap B_2}+...
+a_{B_i\cap B_{C_m^k\cdot C_n^l}}),\eqno(8)$$
and (7) is
$$g_{B_i}=(g_{B_i\cap B_1}\cdot g_{B_i\cap B_2}\cdot ...
\cdot g_{B_i\cap B_{C_m^k\cdot C_n^l}})^{\frac{1}{C_m^k\cdot C_n^l}}.\eqno(9)$$

By the arithmetic-geometric inequality, it follows that
$$a_{B_i\cap B_j} \geq g_{B_i\cap B_j}, ~j=1,2,...,C_m^k\cdot C_n^l.\eqno(10)$$
From (8) and (10), we infer that
$$a_{B_i}\geq\frac{1}{C_m^k\cdot C_n^l}(g_{B_i\cap B_1}+g_{B_i\cap B_2}+...
+g_{B_i\cap B_{C_m^k\cdot C_n^l}}).$$

Therefore
$$\begin{array}{rl}
\Big(\displaystyle\prod_{i=k,j=l\atop B_{ij}\subset B}a_{B_{ij}}\Big)^{\frac{1}{C_m^k\cdot C_n^l}}
&=\Big(\displaystyle\prod_{i=1}^{{C_m^k\cdot C_n^l}}a_{B_i}\Big)^{\frac{1}{C_m^k\cdot C_n^l}}\\
&\displaystyle\geq\frac{1}{C_m^k\cdot C_n^l}\Big(\prod_{i=1}^{C_m^k\cdot C_n^l}(\sum_{j=1}^{C_m^k\cdot C_n^l}g_{B_i\cap
B_j})\Big) ^\frac{1}{C_m^k\cdot C_n^l}.\\
\end{array}\eqno(11)$$

On the other hand, using the discrete case of H\"{o}lder's
inequality in the form
$$\sum_{k=1}^n\Big(\prod_{j=1}^mx_{jk}\Big)^{\frac{1}{m}}
\leq\left(\prod_{j=1}^n\Big(\sum_{k=1}^mx_{jk}\Big)\right)^{\frac{1}{m}},
$$
where $n, m$ are positive integers and $x_{jk}>0(j, k=1, 2,..., m)$, we obtain
$$\Big(\prod_{i=1}^{C_m^k\cdot C_n^l}\Big(\sum_{j=1}^{C_m^k\cdot C_n^l}g_{B_i\cap
B_j}\Big)\Big)^{\frac{1}{C_m^k\cdot C_n^l}}
\geq\sum_{i=1}^{C_m^k\cdot C_n^l}\Big(\prod_{j=1}^{C_m^k\cdot C_n^l}g_{B_i\cap
B_j}\Big)^{\frac{1}{C_m^k\cdot C_n^l}}.\eqno(12)$$ Combining (9), (11) and
(12), ~it follows that
$$\begin{array}{rl}
\Big(\displaystyle\prod_{i=k,j=l\atop B_{ij}\subset B}a_{B_{ij}}\Big)^{\frac{1}{C_m^k\cdot C_n^l}}
&\geq\displaystyle\frac{1}{C_m^k\cdot C_n^l}\displaystyle\sum_{i=1}^{C_m^k\cdot C_n^l}\Big(\displaystyle\prod_{j=1}^{C_m^k\cdot C_n^l}g_{B_i\cap
B_j}\Big)^{\frac{1}{C_m^k\cdot C_n^l}}\\
&=\displaystyle\frac{1}{C_m^k\cdot C_n^l}\displaystyle\sum_{i=1}^{C_m^k\cdot C_n^l}g_{B_i}=\displaystyle\frac{1}{C_m^k\cdot C_n^l}\Big(\sum_{i=k,j=l\atop
B_{ij}\subset B}g_{B_{ij}}\Big)\\
\end{array}$$
which is just the inequality
(2).~~~~~~~~~~~~~~~~~~~~~~~~~~~~~~~~~~~~~~~~~~~~~~~~~~~~~~~~~~~~~~~~~~~~~~~~~~$\square$

\vskip 10pt Lin Si 

College of Science, Beijing Forestry University,

Beijing, 100083, P.R.China

E-mail: silin@bjfu.edu.cn, lin.si@hotmail.com

\vskip 10pt Suyun Zhao

College of Mathematics \& Computer Science, HeBei University,

Baoding, Hebei 071002, P.R.China

E-mail: zhaosy@mail.hbu.edu.cn

\end{document}